\documentclass{article}

\newtheorem{thm}{Theorem}[section]

\newtheorem{lem}[thm]{Lemma}

\begin{document}

\title{{\bf A New Error Bound for Shifted Surface Spline Interpolation}}         % Enter your title between curly braces
\author{{\bf Lin-Tian Luh} \\Department of Mathematics, Providence University\\ Shalu, Taichung\\ email:ltluh@pu.edu.tw \\ phone:(04)26328001 ext. 15126 \\ fax:(04)26324653}        % Enter your name between curly braces
\date{\today}          % Enter your date or \today between curly braces
\maketitle
{\bf Abstract.} Shifted surface spline is a frequently used radial function for scattered data interpolation. The most frequently used error bounds for this radial function are the one raised by Wu and Schaback in \cite{WS} and the one raised by Madych and Nelson in \cite{MN2}. Both are $O(d^{l})$ as $d\rightarrow 0$, where $l$ is a positive integer and $d$ is the well-known fill-distance which roughly speaking measures the spacing of the data points. Recently Luh raised an exponential-type error bound with convergence rate $O(\omega ^{\frac{1}{d}})$ as $d\rightarrow 0$ where $0<\omega <1$ is a fixed constant which can be accurately computed. Although the exponential-type error bound converges much faster than the polynomial-type error bound, the constant $\omega$ is intensely influenced by the dimension $n$ in the sense $\omega \rightarrow 1$ rapidly as $n\rightarrow \infty$. Here the variable $x$ of both the interpolated and interpolating functions lies in $R^{n}$. In this paper we present an error bound which is $O(\sqrt{d}\omega'^{\frac{1}{d}})$ where $0<\omega'<1$ is a fixed constant for any fixed $n$, and is only mildly influenced by $n$. In other words, $\omega'\rightarrow 1$ very slowly as $n\rightarrow \infty$, and $\omega'<<\omega$, especially for high dimensions. Moreover, $\omega'$ can be accurately computed without slight difficulty.\\
\\
{\bf keywords}:radial basis function, shifted surface spline, error bound.
\section{Introduction}       % Enter section title between curly braces

In the theory of radial basis functions, it's well known that any conditionally positive definite radial function can form an interpolant for any set of scattered data. We make a simple sketch of this process as follows.

Suppose $h$ is a continuous function on $R^{n}$ which is strictly conditionally positive definite of order $m$. For any set of data points $(x_{j},f_{j}),j=1,\ldots,N$, where $X=\{ x_{1},\ldots ,x_{N}\}$ is a subset of $R^{n}$ and the $f_{j}'s$ are real or complex numbers, there is a unique function of the form
\begin{equation}
  s(x)=p(x)+\sum_{j=1}^{N}c_{j}h(x-x_{j})
\end{equation}
, where $p(x)$ is a polynomial in $P_{m-1}^{n}$, satisfying
\begin{equation}
  \sum_{j=1}^{N}c_{j}q(x_{j})=0
\end{equation}
for all polynomials $q$ in $P_{m-1}^{n}$ and 
\begin{equation}
  p(x_{i})+\sum_{j=1}^{N}c_{j}h(x_{i}-x_{j})=f_{i},\ i=1,\ldots , N
\end{equation}
if $X$ is a determining set for $P_{m-1}^{n}$.

A complete treatment of this topic can be seen in \cite{MN1} and many other papers.

The function $s(x)$ is called the $h$-spline interpolant of the data points and is of central importance in the theory of radial basis functions. In this paper $h$ always denotes a radial function in the sense that the value of $h(x)$ is completely determined by the norm $| x|$ of $x$. Here, $P_{m-1}^{n}$ denotes the class of those n-variable polynomials of degree not more than $m-1$.

In this paper we are mainly interested in a radial function called shifted surface spline defined by
\begin{eqnarray}
  h(x)& := &(-1)^{m}(|x|^{2}+c^{2})^{\frac{\lambda}{2}}log(|x|^{2}+c^{2})^{\frac{1}{2}},\ \lambda \in Z_{+},\ m=1+\frac{\lambda}{2},\ c>0, \nonumber \\
      &    & x\in R^{n},\ \lambda,n\ even,
\end{eqnarray} 
where $|x|$ is the Euclidean norm of $x$, and $\lambda,c$ are constants. In fact, the definition of shifted surface spline covers odd dimensions. For odd dimensions, it's of the form
\begin{eqnarray}
  h(x)& := & (-1)^{\lceil \lambda -\frac{n}{2} \rceil }(|x|^{2}+c^{2})^{\lambda -\frac{n}{2}},\ n\ odd,\ \lambda \in Z_{+}=\{ 1,2,3,\ldots \} \nonumber \\
      &    & and\ \lambda >\frac{n}{2}.
\end{eqnarray}
However, this is just multiquadric and we treat it in another paper \cite{Lu4}. Therefore we will not discuss it. Instead, we will focus on (4) and even dimensions only. 
\subsection{Polynomials and Simplices}    % Enter subsection title between curly braces
Let $E$ denote an n-dimensional simplex\cite{Fl} with vertices $v_{1},\ldots ,v_{n+1}$. If we adopt barycentric coordinates, then any point $x\in E$ can be written as a convex combination of the vertices:
$$x=\sum_{i=1}^{n+1}\lambda_{i}v_{i},\ \sum_{i=1}^{n+1}\lambda_{i}=1,\ \lambda_{i}\geq 0.$$
 We define the ``{\bf equally spaced}'' points of degree k to be those points whose barycentric coordinates are of the form
$$(k_{1}/k,k_{2}/k,\ldots ,k_{n+1}/k),\ k_{i}\ nonnegative\ integers\ and\ k_{1}+\cdots +k_{n+1}=k.$$
It's easily seen that the number of such points is exactly $dimP_{k}^{n}$, i.e., the dimension of $P_{k}^{n}$. In this section we use N to denote $dimP_{k}^{n}$.

The above-defined equally spaced points can induce a polynomial interpolation process as follows. Let $x_{1},\ldots ,x_{N}$ be the equally spaced points in $E$ of degree k. The associated Lagrange polynomials $l_{i}$ of degree $k$ are defined by the condition $l_{i}(x_{j})=\delta_{ij},\ 1\leq i,j\leq N$. For any continuous map $f\in C(E),\ (\Pi_{k}f)(x):= \sum_{i=1}^{N}f(x_{i})l_{i}(x)$ is its interpolating polynomial. If both spaces are equipped with the supremum norm, the mapping
$$\Pi_{k}:C(E)\rightarrow P_{k}^{n}$$
has a well-known norm 
$$\| \Pi_{k}\| =\max _{x}\sum_{i=1}^{N}|l_{i}(x)|$$
which is the maximum value of the famous Lebesgue function. It's easily seen that for any $p\in P_{k}^{n}$,
$$\| p\| _{\infty}:= \max _{x\in E}|p(x)| \leq \| \Pi_{k}\| \max _{1\leq i\leq N}|p(x_{i})|.$$
The next result is important in our construction of the error bound, and we cite it directly from \cite{Bo}.
\begin{lem}
  For the above equally spaced points $\{ x_{1},\ldots ,x_{N}\}$, $\| \Pi_{k}\| \leq \left( \begin{array}{c}
  2k-1 \\ k
\end{array} \right) $. Moreover, as $n\rightarrow \infty, \ \| \Pi_{k}\| \rightarrow \left( \begin{array}{c}
  2k-1 \\ k
\end{array} \right) $.
\end{lem}
Then we need another lemma which must be proven because it plays a crucial role in our development.
\begin{lem}
  Let $Q\subseteq R^{n}$ be an n simplex in $R^{n}$ and $Y$ be the set of equally spaced points of degree $k$ in $Q$. Then, for any point $x$ in $Q$, there is a measure $\sigma$ supported on $Y$ such that 
$$\int p(y)d\sigma(y)=p(x)$$
for all $p$ in $P_{k}^{n}$, and 
$$\int d|\sigma |(y)\leq \left( \begin{array}{c}
                                2k-1 \\ k
                              \end{array} \right).$$
\end{lem}
{\bf Proof}. Let $Y=\{ y_{1},\ldots , y_{N}\} $ be the set of equally spaced points of degree $k$ in $Q$. Denote $P_{k}^{n}$ by $V$. For any $x\in Q$, let $\delta_{x}$ be the point-evaluation functional. Define $T:V\rightarrow T(V)\subseteq R^{N}$ by $T(v)=(\delta_{y_{i}}(v))_{y_{i}\in Y}$. Then $T$ is injective. Define $\tilde{\psi}$ on $T(V)$ by $\tilde{\psi}(w)=\delta_{x}(T^{-1}w)$. By the Hahn-Banach theorem, $\tilde{\psi}$ has a norm-preserving extension $\tilde{\psi}_{ext}$ to $R^{N}$. By the Riesz representation theorem, each linear functional on $R^{N}$ can be represented by the inner product with a fixed vector. Thus, there exists $z\in R^{N}$ with 
$$\tilde{\psi}_{ext}(w)=\sum_{j=1}^{N}z_{j}w_{j}$$
and $\| z\| _{(R^{N})^{*}}=\| \tilde{\psi}_{ext}\| $. If we adopt the $l_{\infty}$-norm on $R^{N}$, the dual norm will be the $l_{1}$-norm. Thus $\| z\| _{(R^{N}))^{*}}=\| z\| _{1}=\| \tilde{\psi}_{ext}\| =\| \tilde{\psi}\| =\| \delta_{x}T^{-1}\|$.

Now, for any $p\in V$, by setting $w=T(p)$, we have 
$$\delta_{x}(p)=\delta_{x}(T^{-1}w)=\tilde{\psi}(w)=\tilde{\psi}_{ext}(w)=\sum_{j=1}^{N}z_{j}w_{j}=\sum_{j=1}^{N}z_{j}\delta_{y_{j}}(p).$$
This gives 
\begin{equation}
  p(x)=\sum_{j=1}^{N}z_{j}p(y_{j})
\end{equation}
where $|z_{1}|+\cdots +|z_{N}|=\| \delta_{x}T^{-1}\|$.

Note that
\begin{eqnarray*}
  \| \delta_{x}T^{-1}\| & = & \sup _{\begin{array}{c}
                                       w\in T(V) \\ w\neq 0
                                     \end{array} }\frac{\| \delta_{x}T^{-1}(w)\| }{\| w\| _{R^{N}}}\\
                        & = & \sup _{\begin{array}{c}
                                       w\in T(V) \\ w\neq 0
                                     \end{array}}\frac{|\delta_{x}p|}{\| T(p)\| _{R^{N}}}\\
                        & \leq & \sup _{\begin{array}{c}
                                          p\in V \\ p\neq 0
                                        \end{array}}\frac{|p(x)|}{\max _{j=1,\ldots ,N}|p(y_{j})|}\\
                        & \leq & \sup _{\begin{array}{c}
                                          p\in V \\ p\neq 0
                                        \end{array}}\frac{\| \Pi_{k}\| \max _{j=1,\ldots ,N}|p(y_{j})|}{\max _{j=1,\ldots ,N}|p(y_{j})|}\\
                        & = & \| \Pi_{k}\| \\
                        & \leq & \left( \begin{array}{c}
                                          2k-1 \\ k
                                        \end{array}      \right) .
                                                                                                  \end{eqnarray*}
Therefore $|z_{1}|+\cdots +|z_{N}|\leq \left( \begin{array}{c}
                                               2k-1 \\ k
                                             \end{array} \right)$ and our lemma follows immediately from (6) by letting $\sigma(\{ y_{j}\} )=z_{j},\ j=1,\ldots ,N$. \hspace{10.6cm}  $\sharp$
\subsection{Radial Functions and Borel Measures}
Our theory is based on a fundamental fact that any continuous conditionally positive definite radial function corresponds to a unique positive Borel measure. Before discussing this property in detail, we first clarify some symbols and definitions. In this paper ${\cal D}$ denotes the space of all compactly supported and infinitely differentiable complex-valued functions on $R^{n}$. For each function $\phi$ in ${\cal D}$, its Fourier transform is
$$\hat{\phi}(\xi)= \int e^{-i<x,\xi>}\phi(x)dx.$$
Then we have the following lemma which is introduced in \cite{GV} but modified by Madych and Nelson in \cite{MN2}.
\begin{lem}
  For any continuous conditionally positive definite function $h$ on $R^{n}$ of order $m$, there are a unique positive Borel measure $\mu$ on $R^{n}\sim \{ 0\}$ and constants $a_{r},\ |r|=2m$ such that for all $\psi\in {\cal D}$,
\begin{eqnarray}
  \int h(x)\psi(x)dx & = & \int \{ \hat{\psi}(\xi)-\hat{\chi}(\xi)\sum_{|r|<2m}D^{r}\hat{\psi}(0)\frac{\xi^{r}}{r!}\}  d\mu(\xi) \nonumber \\
                    &   & +\sum_{|r|\leq 2m}D^{r}\hat{\psi}(0)\frac{a_{r}}{r!}
\end{eqnarray}
, where for every choice of complex numbers $c_{\alpha},\ |\alpha|=m$,
$$\sum_{|\alpha|=m}\sum_{|\beta|=m}a_{\alpha+\beta}c_{\alpha}\overline{c_{\beta}}\geq 0.$$
Here $\chi$ is a function in ${\cal D}$ such that $1-\hat{\chi}(\xi)$ has a zero of order $2m+1$ at $\xi=0$; both of the integrals
$$\int_{0<|\xi|<1}|\xi|^{2m}d\mu(\xi),\ \int_{|\xi|\geq 1}d\mu(\xi)$$
are finite. The choice of $\chi$ affects the value of the coefficients $a_{r}$ for $|r|<2m$.
\end{lem}
\section{Main Result}
In order to show our main result, we need some lemmas, including the famous Stirling's formula.\\
\\
{\bf Stirling's Formula}: $n!\sim \sqrt{2\pi n}(\frac{n}{e})^{n}$.\\
\\
The approximation is very reliable even for small $n$. For example, when $n=10$, the relative error is only $0.83\%$. The larger $n$ is, the better the approximation is. For further details, we refer the reader to \cite{GG} and \cite{GKP}.

\begin{lem}
  For any positive integer k,
$$\frac{\sqrt{(2k)!}}{k!}\leq 2^{k}.$$
\end{lem}
{\bf Proof}. This inequality holds for $k=1$ obviously. We proceed by induction.
\begin{eqnarray*}
  \frac{\sqrt{[2(k+1)]!}}{(k+1)!} & = & \frac{\sqrt{(2k+2)!}}{k!(k+1)}=\frac{\sqrt{(2k)!}}{k!}\cdot \frac{\sqrt{(2k+2)(2k+1)}}{k+1} \\
                                  & \leq & \frac{\sqrt{(2k)!}}{k!}\cdot \frac{\sqrt{(2k+2)^{2}}}{k+1}\leq 2^{k}\cdot \frac{(2k+2)}{k+1}=2^{k+1}. \hspace{4cm}\ \ \ \ \ \ \ \ \  \sharp 
\end{eqnarray*} 
% Set the ending of a LaTeX document

Now recall that the function $h$ defined in (4) is conditionally positive definite of order $m=1+\frac{\lambda}{2}$. This can be found in \cite{Dy} and many relevant papers. Its Fourier transform \cite{GS} is 
\begin{equation}
  \hat{h}(\theta)=l(\lambda,n)|\theta|^{-\lambda-n}\tilde{{\cal K}}_{\frac{n+\lambda}{2}}(c|\theta|)
\end{equation} 
where $l(\lambda,n)>0$ is a constant depending on $\lambda$ and $n$, and $\tilde{{\cal K}}_{\nu}(t)=t^{\nu}{\cal K}_{\nu}(t)$, ${\cal K}_{\nu}(t)$ being the modified Bessel function of the second kind\cite{AS}. Then we have the following lemma.
\begin{lem}
  Let $h$ be as in (4) and $m$ be its order of conditional positive definiteness. There exists a positive constant $\rho$ such that 
\begin{equation}
  \int_{R^{n}}|\xi|^{k}d\mu(\xi)\leq l(\lambda,n)\cdot \sqrt{\frac{\pi}{2}}\cdot n\cdot \alpha_{n}\cdot c^{\lambda-k}\cdot \Delta_{0}\cdot \rho^{k}\cdot  k!
\end{equation}
for all integer $k\geq 2m+2$ where $\mu$ is defined in (6), $\alpha_{n}$ denotes the volume of the unit ball in $R^{n}$, $c$ is as in (4), and $\Delta_{0}$ is a positive constant.
\end{lem}
{\bf Proof}. We first transform the integral of the left-hand side of the inequality into a simpler form.
\begin{eqnarray}
  &      & \int_{R^{n}}|\xi|^{k}d\mu(\xi)  \nonumber \\
  & =    & \int_{R^{n}}|\xi|^{k}l(\lambda,n)\tilde{{\cal K}}_{\frac{n+\lambda}{2}}(c|\xi|)|\xi|^{-\lambda-n}d\xi \nonumber \\
  & =    & l(\lambda,n)c^{\frac{n+\lambda}{2}}\int_{R^{n}}|\xi|^{k-\frac{n+\lambda}{2}}\cdot {\cal K}_{\frac{n+\lambda}{2}}(c|\xi|)d\xi \nonumber \\
  & \sim & l(\lambda,n)c^{\frac{n+\lambda}{2}}\sqrt{\frac{\pi}{2}}\int_{R^{n}}|\xi |^{k-\frac{n+\lambda}{2}}\cdot \frac{1}{\sqrt{c|\xi|}\cdot e^{c|\xi|}}d\xi \nonumber \\
  & =    & l(\lambda,n)c^{\frac{n+\lambda}{2}}\cdot \sqrt{\frac{\pi}{2}}\cdot n\cdot \alpha_{n}\int_{0}^{\infty}r^{k-\frac{n+\lambda}{2}}\cdot \frac{r^{n-1}}{\sqrt{cr}\cdot e^{cr}}dr  \nonumber \\
  & =    & l(\lambda,n)c^{\frac{n+\lambda}{2}}\sqrt{\frac{\pi}{2}}\cdot n\cdot \alpha_{n}\cdot \frac{1}{\sqrt{c}}\int _{0}^{\infty}\frac{r^{k+\frac{n-\lambda -3}{2}}}{e^{cr}}dr  \nonumber \\
  & =    & l(\lambda,n)c^{\frac{n+\lambda}{2}}\sqrt{\frac{\pi}{2}}\cdot n\cdot \alpha_{n}\cdot \frac{1}{\sqrt{c}}\cdot \frac{1}{c^{k+\frac{n-\lambda -1}{2}}}\int_{0}^{\infty}\frac{r^{k+\frac{n-\lambda -3}{2}}}{e^{r}}dr  \nonumber \\
  & =    & l(\lambda,n)\sqrt{\frac{\pi}{2}}\cdot n\cdot \alpha_{n}\cdot c^{\lambda-k}\int _{0}^{\infty}\frac{r^{k'}}{e^{r}}dr \ where\ k'=k+\frac{n-\lambda -3}{2}. \nonumber 
\end{eqnarray}
Note that $k\geq 2m+2=4+\lambda$ implies $k'\geq \frac{n+\lambda+5}{2}>0$.

Now we divide the proof into three cases. Let $k''=\lceil k'\rceil $ which is the smallest integer greater than or equal to $k'$.

Case1. Assume $k''>k$. Let $k''=k+s$. Then 
$$\int_{0}^{\infty}\frac{r^{k'}}{e^{r}}dr\leq \int_{0}^{\infty}\frac{r^{k''}}{e^{r}}dr=k''!=(k+s)(k+s-1)\cdots (k+1)k!$$
and
$$\int_{0}^{\infty}\frac{r^{k'+1}}{e^{r}}dr\leq \int_{0}^{\infty}\frac{r^{k''+1}}{e^{r}}dr=(k''+1)!=(k+s+1)(k+s)\cdots (k+2)(k+1)k!.$$
Note that
$$\frac{(k+s+1)(k+s)\cdots (k+2)}{(k+s)(k+s-1)\cdots (k+1)}=\frac{k+s+1}{k+1}.$$
The condition $k\geq 2m+2$ implies that
$$\frac{k+s+1}{k+1}\leq \frac{2m+3+s}{2m+3}=1+\frac{s}{2m+3}.$$
Let $\rho=1+\frac{s}{2m+3}$. Then
$$\int_{0}^{\infty}\frac{r^{k''+1}}{e^{r}}dr\leq \Delta_{0}\cdot \rho^{k+1}\cdot (k+1)!$$
if $\int_{0}^{\infty}\frac{r^{k''}}{e^{r}}dr\leq \Delta_{0}\cdot \rho^{k}\cdot k!$. The smallest $k''$ is $k_{0}''=2m+2+s$ when $k=2m+2$. Now,
\begin{eqnarray}
  \int_{0}^{\infty}\frac{r^{k_{0}''}}{e^{r}}dr & = & k_{0}''!=(2m+2+s)(2m+1+s)\cdots(2m+3)(2m+2)! \nonumber \\
                                               & = & \frac{(2m+2+s)(2m+1+s)\cdots (2m+3)}{\rho^{2m+2}}\cdot \rho^{2m+2}\cdot (2m+2)! \nonumber \\
                                               & = & \Delta_{0}\cdot \rho^{2m+2}\cdot(2m+2)!   \nonumber \\ 
                                               &   &  where\ \Delta_{0}=\frac{(2m+2+s)(2m+1+s)\cdots (2m+3)}{\rho^{2m+2}}. \nonumber   
\end{eqnarray}
It follows that $\int_{0}^{\infty}\frac{r^{k'}}{e^{r}}dr\leq \Delta_{0}\cdot \rho^{k}\cdot k!$ for all $k\geq 2m+2$.

Case2. Assume $k''<k$. Let $k''=k-s$ where $s>0$. Then 
$$\int_{0}^{\infty}\frac{r^{k'}}{e^{r}}dr\leq \int_{0}^{\infty}\frac{r^{k''}}{e^{r}}dr=k''!=(k-s)!=\frac{1}{k(k-1)\cdots (k-s+1)}\cdot k!$$
and
\begin{eqnarray}
  \int_{0}^{\infty}\frac{r^{k'+1}}{e^{r}}dr & \leq & \int_{0}^{\infty}\frac{r^{k''+1}}{e^{r}}dr \nonumber \\
                                            & = & (k''+1)!=(k-s+1)!=\frac{1}{(k+1)k\cdots (k-s+2)}\cdot (k+1)!. \nonumber  
\end{eqnarray}
Note that
\begin{eqnarray*}
  &      & \left\{ \frac{1}{(k+1)k\cdots (k-s+2)}/\frac{1}{k(k-1)\cdots (k-s+1)}\right\} \\
  & =    & \frac{k(k-1)\cdots (k- s+1)}{(k+1)k\cdots (k-s+2)} \\
  & =    & \frac{(k-s+1)}{k+1} \\ 
  & \leq &1.
\end{eqnarray*}
Let $\rho =1$. Then
$$\int_{0}^{\infty}\frac{r^{k''+1}}{e^{r}}dr\leq \Delta_{0}\cdot \rho^{k+1}\cdot (k+1)!$$
if $\int_{0}^{\infty}\frac{r^{k''}}{e^{r}}dr\leq \Delta_{0}\cdot \rho^{k}\cdot k!$. The smallest $k$ is $k_{0}=2m+2$. Hence the smallest $k''$ is $k_{0}''=k_{0}-s=2m+2-s$. Now,
\begin{eqnarray}
  \int_{0}^{\infty}\frac{r^{k_{0}''}}{e^{r}}dr & = & k_{0}''!=(2m+2-s)!=(k_{0}-s)! \nonumber \\
                                               & = & \frac{1}{k_{0}(k_{0}-1)\cdots (k_{0}-s+1)}\cdot (k_{0}!) \nonumber \\
                                               & = & \Delta_{0}\cdot \rho^{k_{0}}\cdot k_{0}! \ where\ \Delta_{0}=\frac{1}{(2m+2)(2m+1)\cdots (2m-s+3)}. \nonumber 
\end{eqnarray}
It follows that $\int_{0}^{\infty}\frac{r^{k'}}{e^{r}}dr\leq \Delta_{0}\cdot \rho^{k}\cdot k!$ for all $k\geq 2m+2$.

Case3. Assume $k''=k$. Then 
$$\int_{0}^{\infty}\frac{r^{k'}}{e^{r}}dr\leq \int_{0}^{\infty}\frac{r^{k''}}{e^{r}}dr=k!\ \ and\ \ \int_{0}^{\infty}\frac{r^{k'+1}}{e^{r}}dr\leq (k+1)!.$$
Let $\rho=1$. Then $\int_{0}^{\infty}\frac{r^{k'}}{e^{r}}dr\leq \Delta_{0}\cdot \rho^{k}\cdot k!$ for all $k$ where $\Delta_{0}=1$.

The lemma is now an immediate result of the three cases. \hspace{5.7cm} $\sharp$ \\
\\
{\bf Remark}: For the convenience of the reader we should express the constants $\Delta_{0}$ and $\rho$ in a clear form. It's easily shown that\\
(a)$k''>k$ if and only if $n-\lambda>3$,\\
(b)$k''<k$ if and only if $n-\lambda\leq 1$, and\\
(c)$k''=k$ if and only if $1<n-\lambda \leq 3$,\\
where $k''$ and $k$ are as in the proof of the lemma. We thus have the following situations.\\
(a)$n-\lambda>3$. Let $s=\lceil \frac{n-\lambda-3}{2}\rceil $. Then 
$$\rho=1+\frac{s}{2m+3}\ \ and\ \ \Delta_{0}=\frac{(2m+2+s)(2m+1+s)\cdots (2m+3)}{\rho^{2m+2}}.$$
(b)$n-\lambda\leq 1$. Let $s=-\lceil \frac{n-\lambda-3}{2}\rceil $. Then 
$$\rho=1\ \ and\ \ \Delta_{0}=\frac{1}{(2m+2)(2m+1)\cdots (2m-s+3)}.$$
(c)$1<n-\lambda\leq 3$. We have 
$$\rho=1\ \ and \ \ \Delta_{0}=1.$$

Before introducing our main theorem, we must introduce a function space called {\bf native space}, denoted by ${\bf {\cal C}_{h,m}}$, for each conditionally positive definite radial function $h$ of order $m$. If
$${\cal D}_{m}=\{ \phi\in {\cal D}: \ \int x^{\alpha}\phi(x)dx=0\ for\ all\ |\alpha|<m\}$$
, then ${\cal C}_{h.m}$ is the class of those continuous functions $f$ which satisfy
\begin{equation}
  \left| \int f(x)\phi(x)dx\right| \leq c(f)\left\{ \int h(x-y)\phi(x)\overline{\phi(y)}dxdy\right\} ^{1/2}
\end{equation}
for some constant $c(f)$ and all $\phi$ in ${\cal D}_{m}$. If $f\in {\cal C}_{h,m}$, let $\| f\| _{h}$ denote the smallest constant $c(f)$ for which (10) is true. Then $\| \cdot \| _{h}$ is a semi-norm and ${\cal C}_{h,m}$ is a semi-Hilbert space; in the case $m=0$ it is a norm and a Hilbert space respectively. For further details, we refer the reader to \cite{MN1} and \cite{MN2}. This definition of native space is introduced by Madych and Nelson, and characterized by Luh in \cite{Lu1} and \cite{Lu2}. Although there is an equivalent definition \cite{We} which is easier to handle, we still adopt Madych and Nelson's definition to show the author's respect for them.

Now we have come to the main theorem of this paper. 
\begin{thm}
  Let $h$ be as in (4), and $b_{0}$ be any positive number. Let $\Omega$ be any subset of $R^{n}$ satisfying the property that for any $x$ in $\Omega$, there is an $n$ simplex E of diameter equal to $b_{0}$ and $x\in E\subseteq \Omega$. There are positive constants $\delta_{0},c_{1},\omega',\ 0<\omega'<1,$ for which the following is true:If $f\in {\cal C}_{h,m}$, the native space induced by $h$, and $s$ is the h spline that interpolates f on a subset X of $R^{n}$, then
\begin{equation}
  |f(x)-s(x)|\leq c_{1}\sqrt{\delta}(\omega')^{\frac{1}{\delta}}\cdot \| f\| _{h}
\end{equation}
for all x in $\Omega$ and $0<\delta \leq \delta_{0}$ if $\Omega$ satisfies the property that for any x in $\Omega$ and any number r with $\frac{1}{3C}\leq r\leq b_{0}$, there is an n simplex Q with diameter diamQ=r, $x\in Q\subseteq \Omega$, such that for any integer k with $\frac{1}{3C\delta}\leq k\leq \frac{b_{0}}{\delta}$, there is in Q an equally spaced set of centers from X of degree k-1.(Once $\delta$ and k are chosen, Q is in fact an n simplex of diameter $k\delta$ and X can be chosen to consist only of the equally spaced points in Q of degree k-1.) The number C is defined by 
$$C:=\max \left\{ 8\rho',\ \frac{2}{3b_{0}}\right\} ,\ \rho'=\frac{\rho}{c}$$
where $\rho$ and c appear in Lemma2.2 and (4) respectively. Here $\| f\| _{h}$ is the h-norm of f in the native space. The numbers $\delta_{0},c_{1}$ and $\omega'$ are given by $\delta_{0}:=\frac{1}{3C(m+1)}$ where m appears in (4); $c_{1}:=\sqrt{l(\lambda,n)}\cdot (\pi/2)^{1/4}\cdot \sqrt{n\alpha_{n}}\cdot c^{\lambda/2}\cdot \sqrt{\Delta_{0}}\sqrt{3C}\cdot \sqrt{(16\pi)^{-1}}$ where $\lambda$ is as in (4), $l(\lambda,n)$ appears in (8), $\alpha_{n}$ is the volume of the unit ball in $R^{n}$, and $\Delta_{0}$, together with the computation of $\rho$, is defined in Lemma2.2 and the remark following its proof; $\omega':=(\frac{2}{3})^{1/3C}$.
\end{thm}
{\bf Proof}. Let $\delta_{0}$, and $C$ be as in the statement of the theorem. For any $0<\delta\leq \delta_{0}$, we have $0<\delta \leq \frac{1}{3C(m+1)}$ and $0<3C\delta\leq \frac{1}{m+1}$. Since $\frac{1}{m+1}<1$, there exists an integer $k$ such that
$$1\leq 3C\delta k\leq 2.$$
For such $k$, $\delta k\leq \frac{2}{3C}\leq b_{0},\ \frac{1}{3C\delta}\leq k\leq \frac{b_{0}}{\delta}$, and $8\rho'\delta k\leq \frac{2}{3}$.

For any $x\in \Omega$, choose arbitrarily an n simplex $Q$ of diameter $diamQ=\delta k$ with vertices $v_{0},v_{1},\ldots ,v_{n}$ such that $x\in Q\subseteq \Omega$. Let $x_{1},\ldots ,x_{N}$ be equally spaced points of degree k-1 on $Q$ where $N=dimP_{k-1}^{n}$. By (9) and Lemma2.1, whenever $k>m$,
\begin{eqnarray}
  c_{k} & := & \left\{ \int_{R^{n}}\frac{|\xi|^{2k}}{(k!)^{2}}d\mu(\xi)\right\} ^{1/2} \nonumber \\
        & \leq & \sqrt{l(\lambda,n)}\cdot (\pi/2)^{1/4}\cdot \sqrt{n\alpha_{n}}\cdot c^{\lambda/2}\cdot c^{-k}\cdot \sqrt{\Delta_{0}}\cdot (2\rho)^{k}.
\end{eqnarray}
Recall that Theorem4.2 of \cite{MN2} implies that 
\begin{equation}
  |f(x)-s(x)|\leq c_{k}\| f\| _{h}\cdot \int_{R^{n}}|y-x|^{k}d|\sigma|(y)
\end{equation}
whenever $k>m$, and $\sigma$ is any measure supported on $X$ such that 
\begin{equation}
  \int_{R^{n}}p(y)d\sigma(y)=p(x)
\end{equation}
for all polynomials $p$ in $P_{k-1}^{n}$.

Let $\sigma$ be the measure supported on $\{ x_{1},\ldots ,x_{N}\}$ as mentioned in Lemma1.2. We essentially need to bound the quantity
$$I=c_{k}\int_{R^{n}}|y-x|^{k}d|\sigma|(y)$$
only. Thus for $k$ mentioned as above and $\Delta_{0}$ defined in Lemma2.2, by Lemma1.2,
\begin{eqnarray*}
  I & \leq & \sqrt{l(\lambda,n)}(\pi/2)^{1/4}\sqrt{n\alpha_{n}}c^{\lambda/2}c^{-k}\sqrt{\Delta_{0}}(2\rho)^{k}(\delta k)^{k}\left( \begin{array}{c}
                                    2(k-1)-1 \\ k-1
                                  \end{array}\right) \\
    & \sim & \sqrt{l(\lambda,n)}(\pi/2)^{1/4}\sqrt{n\alpha_{n}}c^{\lambda/2}c^{-k}\sqrt{\Delta_{0}}(2\rho)^{k}(\delta k)^{k}\frac{1}{\sqrt{\pi}}\frac{1}{\sqrt{k-1}}2^{2(k-1)}\ by\ Stirling's\ Formula\\
    & \sim & \sqrt{l(\lambda,n)}(\pi/2)^{1/4}\sqrt{n\alpha_{n}}c^{\lambda/2}c^{-k}\sqrt{\Delta_{0}}(2\rho)^{k}(\delta k)^{k}\frac{1}{\sqrt{\pi}}\frac{1}{\sqrt{k}}4^{k-1}\ when\ k\ is \ large\\
    & = & \sqrt{l(\lambda,n)}(\pi/2)^{1/4}\sqrt{n\alpha_{n}}c^{\lambda/2}\sqrt{\Delta_{0}}\frac{1}{\sqrt{\pi}}\frac{1}{\sqrt{k}}\left( \frac{2\rho \delta k}{c}\right) ^{k}\frac{4^{k}}{4}\\
    & = & \sqrt{l(\lambda,n)}(\pi/2)^{1/4}\sqrt{n\alpha_{n}}c^{\lambda/2}\sqrt{\Delta_{0}}\frac{1}{\sqrt{16\pi}}\frac{1}{\sqrt{k}}\left( \frac{8\rho \delta k}{c}\right) ^{k}\\
    & \leq & \sqrt{l(\lambda,n)}(\pi/2)^{1/4}\sqrt{n\alpha_{n}}c^{\lambda/2}\sqrt{\Delta_{0}}\frac{1}{\sqrt{16\pi}}\frac{1}{\sqrt{k}}\left( \frac{2}{3}\right) ^{k}\\
    & \leq & \sqrt{l(\lambda,n)}(\pi/2)^{1/4}\sqrt{n\alpha_{n}}c^{\lambda/2}\sqrt{\Delta_{0}}\frac{1}{\sqrt{16\pi}}\frac{1}{\sqrt{k}}\left( \frac{2}{3}\right) ^{\frac{1}{3C\delta}}\\
    & = & \sqrt{l(\lambda,n)}(\pi/2)^{1/4}\sqrt{n\alpha_{n}}c^{\lambda/2}\sqrt{\Delta_{0}}\frac{1}{\sqrt{16\pi}}\frac{1}{\sqrt{k}}\left[ \left( \frac{2}{3}\right) ^{\frac{1}{3C}}\right] ^{\frac{1}{\delta}}\\
    & = & \sqrt{l(\lambda,n)}(\pi/2)^{1/4}\sqrt{n\alpha_{n}}c^{\lambda/2}\sqrt{\Delta_{0}}\frac{1}{\sqrt{16\pi}}\frac{1}{\sqrt{k}}[\omega']^{\frac{1}{\delta}}\ where\ \omega'=\left( \frac{2}{3}\right) ^{\frac{1}{3C}}\\
    & \leq & \sqrt{l(\lambda,n)}(\pi/2)^{1/4}\sqrt{n\alpha_{n}}c^{\lambda/2}\sqrt{\Delta_{0}}\frac{1}{\sqrt{16\pi}}\sqrt{3C\delta}[\omega']^{\frac{1}{\delta}}\\
    & = & \sqrt{l(\lambda,n)}(\pi/2)^{1/4}\sqrt{n\alpha_{n}}c^{\lambda/2}\sqrt{\Delta_{0}}\frac{\sqrt{3C}}{\sqrt{16\pi}}\sqrt{\delta}[\omega']^{\frac{1}{\delta}}
\end{eqnarray*}
Our theorem thus follows from (13).\hspace{9.5cm}  $\sharp$\\
\\
{\bf Remark}. In the preceding theorem we didn't mention the well-known fill-distance. In fact $\delta$ is in spirit equivalent to the fill-distance $d(Q,X):=\sup_{y\in Q}\min _{1\leq i\leq N}\| y-x_{i}\|$ where $X=\{ x_{1},\ldots ,x_{N}\}$ as mentioned in the proof. Note that $\delta\rightarrow 0$ if and only if $d(Q,X)\rightarrow 0$. However we avoid using fill-distance  because in our approach the data points are not purely scattered. This to some extent seems to be a drawback. However it does not pose any trouble for us both theoretically and practically. The equally spaced centers $x_{1},\ldots ,x_{N}$ in the simplex $Q$ are friendly and easily tractable.
\section{Comparison}
The exponential-type error bound for (4) is presented by Luh in \cite{Lu3} and is of the form
\begin{equation}
  |f(x)-s(x)|\leq c_{1}\omega^{\frac{1}{\delta}}\| f\| _{h}
\end{equation}
where $c_{1}=\sqrt{l(\lambda,n)}(\pi/2)^{1/4}\sqrt{n\alpha_{n}}c^{\lambda/2}\sqrt{\Delta_{0}}$, $\delta$ is equivalent to the fill-distance and 
$$\omega=\left( \frac{2}{3}\right) ^{\frac{1}{3C\gamma_{n}}}$$
where 
$$C=\max \left\{ 2\rho'\sqrt{n}e^{2n\gamma_{n}},\ \frac{2}{3b_{0}}\right\} ,\ \rho'=\frac{\rho}{c}$$
, $\rho$ and $c$ being the same as this paper, $b_{0}$ be the side length of a cube, and $\gamma_{n}$ being defined recursively by 
$$\gamma_{1}=2,\ \gamma_{n}=2n(1+\gamma_{n-1})\ if\ n>1.$$
The constant $c_{1}$ is almost the same as the $c_{1}$ in (11). The number $b_{0}$ plays the same role as the $b_{0}$ of Theorem2.3. However, $\gamma_{n}\rightarrow \infty$ rapidly as $n\rightarrow \infty$. This can be seen by
$$\gamma_{1}=2,\ \gamma_{2}=12,\ \gamma_{3}=78,\ \gamma_{4}=632,\ \gamma_{5}=6330,\cdots $$
The fast growth of $\gamma_{n}$ forces $e^{2n\gamma_{n}}$ and hence $C$ to grow rapidly as dimension $n\rightarrow \infty$. This means that the crucial constant $\omega$ in (15) turns to 1 rapidly as $n\rightarrow \infty$, making the error bound (15) meaningless for high dimensions.

The advantages of our new approach are:first, there is $\sqrt{\delta}$ in (11) which contributes to the convergence rate of the error bound as $\delta\rightarrow 0$; second, the crucial constant $\omega'$ in (11) are only mildly dependent of dimension $n$. Although $\omega'$ dependends on $\rho$ which in turn depends on $n$, the situation is much better. In fact, $\omega'$ can be made completely independent of $n$ by changing $\lambda$ in (4) to keep $n-\lambda \leq 3$. This can be seen in the remark following Lemma2.2. In other words, we have significantly improved the error bound (15).

\end{document}